\documentstyle[12pt,amscd]{amsart}

\theoremstyle{plain} 

\newtheorem{ThNum}{Theorem}
\newtheorem{PropNum}[ThNum]{Proposition}
\newtheorem{LemNum}[ThNum]{Lemma}

\theoremstyle{definition} 

\newtheorem{Example}{Example} 
\newtheorem{Remark}[Example]{Remark}

\newtheorem{Problem}{Problem}

\newcommand{\hyp}{{\bold H}}

\newcommand{\real}{{\bold R}}
\newcommand{\complex}{{\bold C}}
\newcommand{\proj}{{\bold P}}

\newcommand{\integer}{{\bold Z}}

\newcommand{\ball}{{\bold B}} 
\newcommand{\Hom}{\operatorname{Hom}}

\newcommand{\PGL}{\operatorname{PGL}} 
\newcommand{\PSL}{\operatorname{PSL}} 
\newcommand{\SL}{\operatorname{SL}} 
\newcommand{\PU}{\operatorname{PU}}

\newcommand{\Area}{\operatorname{Area}} 
\newcommand{\Aut}{\operatorname{Aut}}

\begin{document}

\title[Complex Hyperbolic Cone Structures on The Configuration Spaces]{%
	Complex Hyperbolic Cone Structures \\ 
	on The Configuration Spaces} 
\author[S. Kojima]{%
	Sadayoshi Kojima} 
\address{%
	Department of Mathematical and Computing Sciences \\
	Tokyo Institute of Technology \\
	Ohokayama, Meguro \\    
	Tokyo 152-8552 Japan} 
\email{%
	sadayosi@@is.titech.ac.jp}  

\date{%
	Version 1.1 on July 2, 1999} 



\begin{abstract}
The space of marked  $n$  distinct points 
on the complex projective line up to 
projective transformations will be called 
a configuration space in this paper.  
There are two families of 
complex hyperbolic structures on the configuration space 
constructed by Deligne-Mostow and by Thurston.  
We first confirm that these families are the same.  
Then in view of the deformation theory 
for real hyperbolic cone 3-manifolds, 
we review the families for small  $n$.  
\end{abstract}
\maketitle


\section{Introduction}

The space of marked  $n$  distinct points 
on the complex projective line  $\complex \proj^1$  up to 
projective transformations will be called 
a configuration space in this paper 
and we denote it by  ${\cal Q}$.  
It admits a structure of a complex manifold of 
dimension  $n-3$, and 
has a long history for attracting many mathematicians.  
We focus in this paper 
only on results related with complex hyperbolic geometry.  

Deligne and Mostow 
construct a family of equivariant maps of the universal cover
of  ${\cal Q}$  to the  $(n-3)$-dimensional 
complex projective space with respect to the action 
of  $\pi_1({\cal Q})$  and the projective transformations  
in \cite{DeligneMostow}.  
It is parameterized by the exponents of an integral 
representation of a  
several variable analogue of the hypergeometric function.  
The main focus of their paper is to 
discuss when the holonomy representation, 
which is shown to lie 
in  $\PU(1, n-3) \subset \PGL_{n-2}(\complex)$  is discrete,  
and to find many complex hyperbolic lattices.  

On the other hand, 
Thurston provides a different construction of 
complex hyperbolic structures on  ${\cal Q}$  
in \cite{Polyhedra}  based 
on euclidean cone structures on  $\complex \proj^1$,  
each of which is assigned to a configuration via 
a generalized Schwarz-Christofell correspondence.  
It is parameterized by the cone angles. 
His approach 
re-discovers complex hyperbolic lattices found 
by Deligne and Mostow. 
Strictly speaking, Thurston constructed 
structures not on  ${\cal Q}$  but rather 
on the quotient of  ${\cal Q}$  by the action 
of remarking cone points with the same cone angles, 
and in fact he found more lattices.  

Although the discovery of lattices 
has been emphasized as a common part of their results, 
they both actually constructed the continuous families of incomplete 
complex hyperbolic structures on  ${\cal Q}$  
which provide lattices in particular cases.  
The first purpose of this paper is to confirm that 
their underlying families of complex hyperbolic structures on  
${\cal Q}$  are the same.  

Deligne and Mostow studied the family 
in view of 
Mumford's compactification in \cite{Mumford}.  
On the other hand, Thurston viewed their completions 
as cone manifolds.  
However, neither papers emphasize deformation theoretic viewpoints. 
Kapovich and Millson pointed out such aspects 
in relation with the study of mechanical linkages in  
\cite{KapovichMillson, KapovichMillson2}.    
The second purpose of this paper 
is to review their 
families as the deformations of complex hyperbolic cone  
structures on  ${\cal Q}$  for small  $n$.    
It is motivated by 
the deformation theory for real hyperbolic cone 
3-manifolds in  
\cite{ThurstonNote, CullerShalen, NeumannZagier, HodgsonKerckhoff, 
Kojima, CooperHodgsonKerckhoff}, 
The study stays still in very primitive stage, 
but a few small, and we believe suggestive, 
observations will be presented.  


\section{Configuration space} 

A configuration of marked  $n$  points on  $\complex \proj^1$  
is the way to distribute points with markings  
on  $\complex \proj^1$  disjointly. 
Let  ${\cal Q}$  be the space of configurations 
of marked  $n$  points up to projective transformations, 
and call it a configuration space.  
That is to say, if we let the space of configurations,   
\begin{equation*} 
	{\cal M} = 
		\overbrace{\complex \proj^1 \times \cdots 
			\times \complex \proj^1}^{%
			\text{$n$ times}} - {\cal D},  
\end{equation*} 
where  ${\cal D}$  is the big diagonal set, then 
\begin{equation*}
	{\cal Q} = {\cal M}/\PGL_2(\complex),  
\end{equation*} 
where  $\PGL_2(\complex)$  acts diagonally.  
By sending the last three marked points to  $\{0, 1, \infty \}$, 
we can always normalize a configuration so that 
the first $n-3$ points lie in  $\complex - \{0, 1 \}$.  
This normalization gives a canonical identification 
of  ${\cal M}$  with the product  
${\cal Q} \times \PGL_2(\complex)$.  
By definition,  ${\cal Q}$  admits a canonical action of 
the symmetry group of  $n$  letters by remarking the points. 

\begin{Example}\label{Ex:Q_4} 
When  $n = 4$, 
${\cal Q}$  is homeomorphic to  $\complex \proj^1 - \{ 0, 1, \infty \}$.  
The action of the symmetry group of markings, 
say $\{1, 2, 3, 4 \}$,  on  ${\cal Q}$  is not 
effective, 
because the action of the Klein permutation group  
$\{ e, (12)(34), (13)(24), (14)(23) \}$  is 
realized by projective transformations.  
The quotient group  $\Gamma$, isomorphic to a dihedral group 
of order  $6$, acts effectively on  ${\cal Q}$.  
${\cal Q}/\Gamma$  is naturally extends to 
an orbifold isomorphic to 
the moduli space  $\hyp/\PSL_2(\integer)$  of 
elliptic curves.  
Such ineffectiveness of the action of the symmetry group occurs only 
when  $n = 4$.  
\end{Example} 

\begin{Example}\label{Ex:Q_5}
Example 1 of  \S4  in \cite{DeligneMostow} discusses 
what  ${\cal Q}$  looks like when  $n = 5$.  
It can be identified with the 
complement of seven rational curves 
in  $\complex \proj^1 \times \complex \proj^1$  
defined below, 
\begin{equation*}
	x = 
		\begin{cases} 
			\quad 0 \\
			\quad 1\\
			\quad \infty, 
		\end{cases} 
	\qquad 
	y = 
		\begin{cases} 
			\quad 0 \\
			\quad 1\\
			\quad \infty, 
		\end{cases} 
	\qquad 
	x = y, 
\end{equation*} 
where  $(x, y) \in \complex \proj^1 \times \complex \proj^1$.  
$(0,0), \, (1,1)$  and  $(\infty,\infty)$  are 
the points where three curves meet, see Figure 1.  
To get a more symmetric representative
with respect to 
the action of the symmetry group of five letters, 
we may blow up these three points.  
Then  ${\cal Q}$  is homeomorphic to 
the complement of ten  $-1$  rational curves in 
$(\complex \proj^1 \times \complex \proj^1) \# 
3 \overline{\complex \proj^2} \approx 
\complex \proj^2 \# 4 \overline{\complex \proj^2}$.  
\begin{figure} 
	\begin{center} 
		\setlength{\unitlength}{1mm} 
		\begin{picture}(80,80) 
		\multiput(10,19)(0,21){3}{\line(1,0){60}}
		\multiput(19,10)(21,0){3}{\line(0,1){60}}
		\put(10,10){\line(1,1){60}}
		\put(19,19){\circle*{2}}
		\put(40,40){\circle*{2}}
		\put(61,61){\circle*{2}}
\end{picture}
	\end{center} 
	\caption{${\cal Q}$ \; for $n=5$}
\end{figure}
\end{Example} 

The complex hyperbolic structure on  ${\cal Q}$  by 
Deligne-Mostow to be discussed depends on the weight 
which will be described by a vector of real numbers, 
\begin{equation}\label{Eq:weight}
	\mu = (\mu_1, \mu_2, \cdots, \mu_n) 
	\quad \text{such that} \quad
	0 < \mu_j < 1 \quad \text{and} \quad  \sum_j \mu_j = 2.  
\end{equation} 
This appears soon as exponents of some multi-valued 1-form.  
It is related with an angle vector 
\begin{equation*}
	\theta = (\theta_1, \theta_2, \cdots, \theta_n) 
	\quad \text{such that} \quad
	0 < \theta_j < 2 \pi \quad \text{and} \quad  
		\sum_j (2 \pi - \theta_j) = 4 \pi   
\end{equation*} 
in Thurston's complex hyperbolization subject to  
the identity, 
\begin{equation*} 
	\theta_j = 2 \pi(1 - \mu_j).  
\end{equation*} 
The weight  $\mu$  can be regarded as a 
curvature vector from Thurston's viewpoint.  

To construct structures in both methods,  
the common root is an integrand of 
an integral representation of a several 
variable analogue of the hypergeometric function  
\begin{equation}\label{Eq:Hypergeometric} 
	\omega_m = \prod (z - m_j)^{-\mu_j} dz    
\end{equation}
assigned to each configuration 
\begin{equation*} 
	m = (m_1, m_2, \cdots, m_n) \in  {\cal M}.   
\end{equation*} 
If one of  $m_j$'s  is  $\infty$, 
we should appropriately understand the 
representation  (\ref{Eq:Hypergeometric})  
as carefully explained in \cite{DeligneMostow}.  
We will see their constructions more precisely in 
the next two sections.


\section{Deligne-Mostow's construction} 

Let  $\proj_m$  be the complement of the point set  
$\{m_1, m_2, \cdots, m_n\}$  in  $\complex \proj^1$, 
namely 
\begin{equation*} 
	\proj_m = \complex \proj^1 - \{m_1, \cdots, m_n\}.  
\end{equation*}   
The construction by Deligne and Mostow in  \cite{DeligneMostow}   
starts with 
choosing a flat complex line bundle  $L_m$  on  $\proj_m$  
with holonomy so that 
the image of a tiny circle surrounding the point 
marked by  $m_j$  
is the rotation of  $2 \pi \mu_j$.  
In other words, 
the holonomy around  $m_j$  acts on the fiber as 
a complex multiplication by  $e^{2 \pi i \mu_j}$. 
$L_m$  admits a hermitian structure, and  
we choose one, though the structure is not unique since 
$\Aut L_m$  is isomorphic to  $\complex^*$.  
The monodromy of  $\omega_m$  around  $m_j$  
is the inverse of that of a horizontal section of  
$L_m$.  
Hence any section of  $\Omega^1(L_m)$  can 
be written as a tensor product of  $\omega_m$, 
a non zero multi-valued section of  $L_m$  and 
a holomorphic function on  $\proj_m$.  

Then consider de Rham cohomology of 
$\proj_m$  with coefficients in  $L_m$.  
Since  $L_m$  is nontrivial by definition of  $\mu$, 
the zero-th cohomology vanishes.  
Thus by Euler characteristic argument, 
the first cohomology group is an  $(n-2)$-dimensional 
complex vector space.  
The hermitian structure we put on  $L_m$  defines 
a hermitian structure on  $H^1(\proj_m; L_m)$.  

Since each  $\mu_j$  lies 
between  $0$  and  $1$,  or the rotation angles lie between  
$0$  and  $2\pi$,  
Proposition 2.6.1 in \cite{DeligneMostow} identifies 
the cohomology group in question with 
that with compact support by the induced 
homomorphism of the inclusion.  
Namely 
\begin{equation*}
	H^1(\proj_m; L_m) \cong 
		H^1_c(\proj_m; L_m) 
\end{equation*}

Poincar\'e duality pairing in this setting defines 
a perfect pairing 
\begin{equation*} 
	\psi_0 : H_c^1(\proj_m;L_m) \times 
		H_c^1(\proj_m;\overline{L_m}) \to 
			H_c^2(\proj_m ;\complex) \cong \complex 
\end{equation*} 
by sending  $\omega_1 \in H_c^1(\proj_m;L_m)$  and 
$\omega_2 \in H_c^1(\proj_m;\overline{L_m})$  to 
\begin{equation*} 
	\psi_0 (\omega_1, \omega_2) = 
		\int_{\proj_m} \omega_1 \wedge \omega_2,    
\end{equation*} 
where  $\overline{L_m}$  is the complex conjugate to  $L_m$.   
This now gives a hermitian form
\begin{equation*} 
	\psi(\omega, \eta) =
		\frac{-1}{2 \pi i} \psi_0(\omega, \overline{\eta})  
\end{equation*} 
on  $H^1_c(\proj_m;L_m)$.  
Corollary 2.21 in \cite{DeligneMostow} shows that  
the hermitian form  $\psi$  is nondegenerate and 
has signature  $(1, n-3)$  
by the Hodge theory.  
Moreover  $\omega_m$  represents a non zero class 
which lies in the positive part 
with respect to  $\psi$  in  $H^1_c(\proj_m ; L_m)$.  

Let  $U$  be a contractible neighborhood of  $m$  in  
${\cal M}$.  
Then  $L_m$  extends uniquely to a flat line bundle  
$L_U$  on  $\cup_{m \in U} \proj_m$.  
Regarding  $L_U$  as a sheaf of horizontal sections, 
and taking a higher direct image of  $L_U$  of 
the projection  $\pi : \cup_{m \in U} \proj_m \to U$, 
we obtain a sheaf  $R^1 \pi_* L_U$  on  $U \subset {\cal M}$  
whose stalk at  $m$  is identified with 
a vector space  $H_c^1(\proj_m ; L_m)$.  
Hence  $R^1 \pi_* L_U$  can be viewed also as 
a flat vector bundle.  
Now the flat projective space bundle  
$P R^1 \pi_* L_U$  is independent 
of the choice of  $L_U$  up to unique isomorphism, 
and hence for variable   $U$,  
they glue into a flat projective space bundle 
on the whole  ${\cal M}$. 
We denote this flat projective space bundle 
by  $B(\mu)$  where the fiber is the projective space 
of the first cohomology group  $H^1_c$.  

Lemma 3.5 in  \cite{DeligneMostow}  shows that 
the assignment of  $[\omega_m]$  to 
each  $m \in {\cal M}$  defines 
a holomorphic section  
\begin{equation*} 
	\omega_{\mu} : {\cal M} \to B(\mu) 
\end{equation*} 
which is equivariant with respect to 
the action of  $\PGL_2(\complex)$.  
Hence restricting  $\omega_{\mu}$  to 
${\cal Q}$,  
we get a section on  ${\cal Q}$  
\begin{equation*} 
	\omega_{\mu} \vert_{{\cal Q}} : 
		{\cal Q} \to B(\mu) \vert_{{\cal Q}} 
\end{equation*} 

Let  $p : \widetilde{\cal Q} \to {\cal Q}$  be the 
universal covering.  
Then the pull back  $p^* B(\mu) \vert_{\cal Q}$  admits 
the product structure   $\widetilde{\cal Q} \times B(\mu) \vert_0$  
induced by the flat structure, 
where  $0$  denotes a fixed base configuration 
lying in  ${\cal Q} \subset {\cal M}$.  
Hence composing the pull back of  $\omega_{\mu}$  and 
the projection  $: \widetilde{\cal Q} \times B(\mu) \vert_0 
\to B(\mu) \vert_0$,  
we get a map 
\begin{equation*} 
	\widetilde{\omega_{\mu}} : \widetilde{\cal Q} 
		\to B(\mu) \vert_0 = \complex \proj^{n-3}.      
\end{equation*} 
Proposition 3.9 in \cite{DeligneMostow} establishes that 
$\widetilde{\omega_{\mu}}$  is locally biholomorphic.  
Moreover (3.10) in \cite{DeligneMostow} shows that 
the image of  $\widetilde{\omega_{\mu}}$  is contained in 
the complex ball  $\ball \subset B(\mu) \vert_0$,  where    
$\ball$  is the quotient of positive part of  $\psi$  
by  $\complex^*$ action.  
Thus the action of  $\pi_1({\cal Q})$  on  $\complex \proj^{n-3}$  
is contained in  $\PU(1,n-3)$  and 
$\widetilde{\omega_{\mu}}$  is 
equivariant with respect to the action of  $\pi_1({\cal Q})$

To end the construction, 
notice that  $\psi$  induces a Bergman metric on  $\ball$  which 
we call a complex hyperbolic metric.  
Pull back this metric on  $\widetilde{\cal Q}$  by 
$\widetilde{\omega_{\mu}}$.   
Since the holonomy representation of  $\pi_1({\cal Q})$  
preserves the metric,  
the metric on  $\widetilde{\cal Q}$  is preserved by 
the action of the covering transformations.  
Hence it descends to a complex hyperbolic structure 
on  ${\cal Q}$.  
The structure depends continuously on  $\mu$, 
and hence we obtained a family of complex hyperbolic   
structures on  ${\cal Q}$  parameterized by the 
weight  $\mu$.  
This summarizes the construction by Deligne and Mostow.  

Fixing  $\mu$, we thus obtained a complex hyperbolic 
structure on  ${\cal Q}$.  
Let us denote by  ${\cal MD}(\mu)$  the completion of 
a complex hyperbolic manifold so constructed.


\section{Thurston's construction}

The method of complex hyperbolization by 
Deligne and Mostow involves the complex Lorentz 
space supported on the first cohomology group 
of  $\proj_m$  with 
a twisted coefficient  $L_m$  together with 
a hermitian form  $\psi$  derived from Poincar\'e duality pairing.  
Thurston gave a completely different aspect of 
these machineries.  
Here we describe how he translated these ideas to his own. 

Fixing a base point  $\ast$  in  $\proj_m$, 
Thurston regards the integral of  $\omega_m$  
along a path from  $\ast$  to  $z$  in  $\proj_m$,     
\begin{equation*} 
	h(z) = \int_{\ast}^z \, \omega_m  
		= \int_{\ast}^z \, \prod (t - m_j)^{-\mu_j} dt  
\end{equation*}
as a developing map of some euclidean structure 
on  $\proj_m$  which extends to an euclidean cone 
structure on  $\complex \proj^1$  with prescribed cone data, 
and relate the family of euclidean cone spheres obtained by 
varying  $m$  with a complex hyperbolic structure on  ${\cal Q}$.  

The reason why the euclidean cone structure 
appears comes from the fact that  
the pre Schwarzian of a multi-valued map  $h$  has 
the form 
\begin{equation*} 
	\frac{h''}{h'} = \sum_j \frac{-\mu_j}{z-m_j}, 
\end{equation*} 
and is single-valued.  
This fact implies that the change of the analytic continuation around 
singular point  $m_j$  is a post composition of 
a map which is necessarily affine.  
Moreover direct computation shows that 
the map must preserve an euclidean metric.  

Proposition 6.1 in \cite{Polyhedra} shows a method to 
assign to each configuration an euclidean cone sphere as follows. 
Fix an euclidean metric on  $\complex$.  
For each configuration  $m \in {\cal M}$, 
we choose a representative such that 
non of  $m_j$'s  is  $\infty$.  
Computation shows that the pre Schwarzian 
is holomorphic at  $\infty$  since  
$\sum_j \mu_j = 2$.  
Thus  $h$  defines a $\pi_1(\proj_m)$-equivariant 
map of the universal cover of  $\proj_m$  to 
$\complex$.  
The image of the holonomy representation 
is contained in the group of euclidean isometries.  
By pulling back the euclidean metric of  $\complex$  on 
the universal cover of  $\proj_m$,  and pushing down to  $\proj_m$,  
we get an euclidean metric there.  
The metric is not complete, 
and the completion 
yields a cone point of cone angle  $2\pi(1-\mu_j)$,  
or curvature  $\mu_j$,  
at each punctured point.  
We denote such an  
euclidean cone sphere by  $\Delta_m$.  

This correspondence is not quite one to one since 
there are several choices we made.  
However, 
it turns out to be one to one 
if we regard it as a correspondence between 
the set of projective classes of configurations with weight  $\mu$   
and the set of similarity classes of euclidean 
cone spheres with prescribed 
curvature  $\mu$.   
In fact, 
the converse is obtained by 
remembering only a conformal structure on 
$\proj_m$  induced from an euclidean structure 
and extend it to the unique conformal structure on $\complex \proj^1$.  

The method of complex hyperbolization by Thurston is 
to give a local coordinate around  $\Delta_m$. 
To do this, choose a geodesic triangulation  $T$  of
$\Delta_m$  such that 
vertices consists of cone points.  
Such a triangulation certainly exists 
by Proposition 2.1 in \cite{Polyhedra}.  
Fixing a triangulation  $T$,  we consider the set  $E$  of 
oriented edges of the universal cover of  
$\Delta_m - \text{\{ cone points \}}$.  
Assigning to each edge in  $E$  the difference of 
the images of the end point and the terminal point by  $h$, 
we get a map  $z_m : E \to \complex$.  
The map  $z_m$  satisfies the following 
cocycle conditions with twisted coefficients in  $L_m$, 
\begin{enumerate} 
	\item 
	$z_m(e_1) + z_m(e_2) + z_m(e_3) = 0$,  \quad 
	when $e_1, e_2, e_3$ surround a triangle, 
	\item 
	$z_m(\gamma e) = H(\gamma) z_m(e)$,  \quad 
	where  $H(\gamma)$  is a rotation part of 
	the holonomy of $\gamma$
\end{enumerate} 
This is well defined up to  $\complex^*$  action.  
Note that the rotation part  $H$  depends only on  
the curvature  $\mu$  and not on the location of 
cone points  $m$.  

The set of euclidean cone spheres 
close to  $\Delta_m$  up to similarity can be parameterized 
locally by cocycles such as 
\begin{equation*} 
	Z = \{ z : E \to \complex \, \vert \, 
		z(e_1)+z(e_2)+z(e_3)=0, 
		\; z(\gamma e) = H(\gamma) z(e) \} 
\end{equation*} 
up to  $\complex^*$ action.  
Proposition 2.2 in \cite{Polyhedra} shows that  $Z$  is 
a complex vector space 
of dimension  $n-2$,  and each cocycle can 
be determined by choosing the values of 
$n-2$ edges  $e_1, e_2, \cdots, e_{n-2}$  
which form a tree in  $E$  and also in  $\Delta_m$.  

\begin{LemNum}\label{Lem:Bijection}
The assignment of  $\omega_m \in B(\mu) \vert_0$  to 
$z_m \in PZ$  provides a local bijection, 
where  $PZ$  is a projective space of  $Z$.  
\end{LemNum} 

\begin{pf}  
Fix a configuration  $m_0$.  
Then  $z_m$  near  $z_{m_0}$  is parameterized by the value of 
appropriate  $n-2$  edges  $e_1, \cdots, e_{n-2}$  up to 
$\complex^*$ action and hence  
$(z(e_1), \cdots, z(e_{n-2}))$  provides its 
virtual coordinate in  $\complex^{n-2}$.  
Easy calculation shows 
\begin{equation*} 
	z_m(e_j) = \int_{h(e_j)} \, dz  
		= \int_{e_j} h^* \, dz 
		= \int_{e_j} h' \, dz 
		= \int_{e_j} \omega_m.   
\end{equation*} 
On the other hand, 
identifying  $\Delta_m$  with a conformal 
extension of  $\proj_m$  to  $\complex \proj^1$, and 
listing the evaluation of  $\omega_m$  along the edges  
$e_1, \cdots, e_{n-2}$  in the last term,  
we get a period integral.  
which induces a virtual coordinate 
of  $\omega_m$  in  $B(\mu)\vert_0$.  
\end{pf} 

Assigning the area of  $\Delta_m$  
to each cocycle  $z_m$, 
we get a hermitian form  $\Area$  on  $Z$, 
\begin{equation*} 
	\Area : Z \to \real \subset \complex. 
\end{equation*}  
Proposition 2.3 in \cite{Polyhedra}  shows that  
$\Area$  turns out to be a hermitian form of signature  $(1, n-3)$, 
and hence induces a complex hyperbolic metric on the ball in  $PZ$.  

Each cocycle under the triangulation gives a virtual local chart 
up to  $\complex^*$ action.  
The coordinate change is attained by changing triangulations.  
However  $\Area$  is invariant under the coordinate change 
up to  $\complex^*$ action.  
hence the system of coordinate charts so constructed 
defines a complex hyperbolic structure on  ${\cal Q}$.  
We denote its completion by  ${\cal T}(\mu)$.  

\begin{LemNum} 
$\Area$  equals  $\pi \psi$  by 
the correspondence in Lemma \ref{Lem:Bijection}.  
\end{LemNum} 

\begin{pf} 
It is enough to verify the identity for 
a geodesic triangle $\Delta$  on  $\Delta_m$.  
The area of  $\Delta$  is equal by definition to 
\begin{equation*} 
	\frac{-1}{2i} \int_{h(\Delta)} dz \wedge d \overline{z} 
	= \frac{-1}{2i} \int_{\Delta} h^* (dz \wedge d \overline{z}) 
	= \frac{-1}{2i} \int_{\Delta} \vert h'(z) \vert^2 
		dz \wedge d \overline{z} 
	= \frac{-1}{2i} \int_{\Delta} \omega_m \wedge \overline{\omega_m}.  
\end{equation*} 
\end{pf} 

\begin{ThNum}\label{Th:DM=T} 
${\cal DM}(\mu)$  is canonically 
isometric to  ${\cal T}(\mu)$.  
\end{ThNum} 

\begin{pf} 
Fix the weight or curvature  $\mu$.  
Then since the local charts of Deligne-Mostow and Thurston 
for  $\cal Q$  are equivalent, 
and the metrics they put are the same, 
they are isometric.  
So are their completions.     
\end{pf}

\begin{Remark} 
As mentioned in the introduction, 
Thurston constructed a complex hyperbolic 
structure not on  ${\cal Q}$  but on 
the quotient of  ${\cal Q}$  by 
the action of remarking 
the cone points with the same cone angles.  
Hence very precisely speaking, 
${\cal T}(\mu)$  agrees with his only when 
cone angles all are mutually distinct. 
\end{Remark}


\section{Deformations} 

Both constructions provide a family of incomplete 
complex hyperbolic structures on  ${\cal Q}$.  
Deligne-Mostow discussed the compactification 
in relation with 
Mumford's geometric invariant theory \cite{Mumford}.  
In particular, topological stratification of the completion 
has been clarified. 
For example, the role of stable and semistable points 
is extensively studied in 
\S\S 6-7 in \cite{DeligneMostow}.  
On the other hand, 
Thurston discussed the completion from 
geometric viewpoints by introducing 
complex hyperbolic cone structures.  
For example, he showed 

\begin{PropNum}[Proposition 2.5 
in  \cite{Polyhedra}]\label{Prop:ConeAngle} 
The cone angle around the complex codimension one singularity 
arisen as collisions of two points with curvature  
$\mu_j, \mu_i$  such that  $\mu_j + \mu_i \leq 1$  
is  $2\pi (1 - \mu_j - \mu_i)$.  
\end{PropNum} 
  
The family provides the deformations of 
complex hyperbolic cone structures on fairly 
stable underlying topological space. 
We will look at them from deformation theoretic viewpoint in 
this section.  

By virtue of Theorem \ref{Th:DM=T}, 
we denote both  
${\cal DM}(\mu)$  and  ${\cal T}(\mu)$   
by  ${\cal Q}(\mu)$.  

Start with the classical case when  $n = 4$. 
Recall that  ${\cal Q}$  is homeomorphic to 
$\complex \proj^1 - \{ 0, 1, \infty \}$.  
${\cal Q}(1/2,1/2,1/2,1/2)$  is 
isometric to a hyperbolic surface homeomorphic to a three 
punctured sphere.  
When the weight varies to 
$\mu = (\mu_1, \mu_2, \mu_3, \mu_4)$,  then 
${\cal Q}(\mu)$  becomes a hyperbolic cone sphere. 
Since the total sum of  $\mu_j$'s  equals  $2$, 
at most three pairs of  $\mu_j$'s  have 
the sum  $\mu_j + \mu_i$  less than  $1$.  
Such a pair provides a cone singularity of cone angle 
$= 2\pi(1-\mu_j-\mu_i)$.  
If there are less than three such pairs, 
then there are pairs whose sum equals  $1$.  
Such a pair provides a cusp.  
The total number of cusps and cone points must be three.  

\begin{ThNum} 
Any real hyperbolic cone sphere with $3$ cone points 
(including cusps) whose cone angles all are less than $2 \pi$  occurs 
as  ${\cal Q}(\mu)$  for some  $\mu =(\mu_1, \mu_2, \mu_3, \mu_4)$.  
\end{ThNum} 

\begin{pf}
The isometry classes of hyperbolic cone spheres 
with three cone points are classified by the 
cone angles.  
Hence it is sufficient to 
solve an equation, for example, 
\begin{equation*} 
	\begin{cases} 
		\; A = 2\pi (1 - \mu_2 - \mu_3), \\
		\; B = 2\pi (1 - \mu_3 - \mu_1), \\
		\; C = 2\pi (1 - \mu_1 - \mu_2),   
	\end{cases} 
\end{equation*}  
for given nonnegative constants  $A, B, C$  
such that  $0 \leq A + B + C < 2 \pi$, and 
to let  $\mu_4 = 4\pi - 2 \pi (\mu_1+\mu_2+\mu_3)$.   
\end{pf}

\begin{Example}\label{Ex:IsometricResult} 
Different weights still 
can give isometric cone spheres in this case.  
For instance, 
${\cal Q}(1/2-\varepsilon, 1/2-\varepsilon, 
1/2-\varepsilon, 1/2+3\varepsilon)$  and 
${\cal Q}(1/2-3\varepsilon, 1/2+\varepsilon, 
1/2+\varepsilon, 1/2+\varepsilon)$  both 
give a hyperbolic cone sphere with three 
cone points of cone angle  $= 2\varepsilon$.   
These weights cannot be transformed by any permutation 
of markings.  
\end{Example} 

When  $n = 5$, the situation is a bit complicated.  
Recall that  ${\cal Q}$  is homeomorphic to 
the complement of the union of ten $-1$  rational curves 
in  $X = (\complex \proj^1 \times \complex \proj^1) {\#} 3 
\overline{\complex \proj^2} \approx 
\complex \proj^2 {\#} 4 \overline{\complex \proj^2}$  as in 
Example \ref{Ex:Q_5}.  
We denote the union of these curves by  ${\cal L}$.  
The pair  $(X, {\cal L})$  will be 
a basic underlying topological space of complex hyperbolic manifolds 
we discuss.  
There is a natural way to index each irreducible component of 
${\cal L}$  by  ${\cal L}_{ji}$  where  $j, i$  are 
integers such that  $1 \leq j < i \leq 5$.  
The index has the property that 
${\cal L}_{ji}$  does intersect with  ${\cal L}_{kl}$  
iff  $\{ j, i \} \cap \{ k, l \} = \emptyset$.  
In fact,  ${\cal L}_{ji}$  can be 
identified with the set of degenerate configurations 
by the collision of the points marked by  $m_j$  and  $m_i$  
under some weight  $\mu$.  

\begin{Example}\label{Ex:TopologicallyStableDeformation}
${\cal Q}(2/5,2/5,2/5,2/5,2/5)$  is a 
compact complex hyperbolic cone manifold, 
where the singular set is located exactly as 
${\cal L} = \cup_{ji} {\cal L}_{ji}$.  
The cone angles around  ${\cal L}_{ji}$  all are  $2\pi/5$  and 
hence it is an orbifold.  
The intersection of  ${\cal L}_{ji}$  and  ${\cal L}_{kl}$  
if any corresponds to the simultaneous collision of 
two pair of points. 
When  $\mu = (\mu_1, \cdots, \mu_5)$  varies such that 
\begin{equation}\label{Eq:Stability} 
	\mu_j + \mu_i < 1 \quad \text{for all} \;\;  i \ne j, 
\end{equation} 
then the underlying topology of  ${\cal Q}(\mu)$  is 
stable and the pair with the singular set is homeomorphic to 
$(X, {\cal L})$. 
The cone angle around  ${\cal L}_{ji}$  is 
equal to  $2\pi(1-\mu_j-\mu_i) \, (< 2\pi)$  by 
Proposition \ref{Prop:ConeAngle}  and 
it is easy to see that 
the parameter space of  $\mu$  under the condition  (\ref{Eq:Stability}) 
injects into the space of marked complex hyperbolic cone structures 
on  $(X, {\cal L})$  by looking at cone angles 
appeared in  ${\cal Q}(\mu)$.   
Similar injectivity can be established for 
odd  $n$  under the condition (\ref{Eq:Stability}). 
\end{Example} 

To see the limiting case and beyond when  $n = 5$, 
we briefly review what happens in real dimension  $3$. 
There are essentially two types of corresponding 
deformations in real hyperbolic cone 3-manifolds, 
which are cusp openings.  

One is provided by throwing a geodesic cone singularity 
away to  $\infty$  and opening a cusp, 
which was discussed originally in \cite{ThurstonNote} 
as a part of the hyperbolic Dehn filling theory, and 
studied as a deformation of cone manifolds in  \cite{Kojima}.  
This is due to the existence of codimension two 
euclidean line.  
In this case, 
the continuous deformations beyond the limit may be regarded 
as cone manifolds with different topology,.  
Some particular discussions of such deformations 
related with the configuration space can be found in 
\cite{KojimaNishiYamashita}.  

The other example is discussed in Example 7.2 in \cite{Kojima}.  
It is provided by collapsing a totally 
geodesic hyperbolic cone sphere in 
a real hyperbolic cone 3-manifold 
to a splitting euclidean cone sphere.  
This is due to the existence of geodesic hypersurfaces. 
In this case, 
the continuous deformations beyond the limit may 
be regarded as one having a vertex singularity 
where the cone axis which were stuck through 
the cone sphere meet.  

The complex hyperbolic geometry of  $\dim_{\complex} \geq 2$  does 
not admit neither real geodesic hypersurfaces, nor 
real codimension two euclidean surfaces.  
Hence it is not conceivable to expect 
a direct analogue of a cusp opening 
deformations in the real case.  
However one sees below that we certainly 
have cusp opening deformations when the 
condition (\ref{Eq:Stability})  breaks down.  
It can be understood as a mixed type of 
two cases in real dimension 3.  

\begin{Example}\label{Ex:CuspOpeningInComplexHyperbolicGeometry}
When the weight approaches $(1/2,1/2,1/3,1/3,1/3)$,  
then the cone angle 
around  ${\cal L}_{12}$  becomes zero and   
${\cal L}_{12}$  itself escapes away to the cusp.  
${\cal L}_{12}$  is topologically a sphere in  
${\cal Q}(2/5,2/5,2/5,2/5,2/5)$  and metrically 
a hyperbolic cone sphere with 
three cone points of cone angles  $2\pi/5$.  
Since the first Chern class of the normal 
bundle of  ${\cal L}_{12}$  is  $-1$, 
the boundary of an equidistant neighborhood of 
${\cal L}_{12}$  supports  $\widetilde{\SL_2(\real)}$ geometry.  
According to the deformation,  
${\cal L}_{12}$  approaches  $\infty$  
where its rescaling limit is an euclidean cone sphere with 
three cone points of cone angle  $2\pi/3$,  and 
the section at the cusp supports nilgeometry.  

When the weight goes beyond the point, say 
$\mu = (5/8,5/8,1/4,1/4,1/4)$, 
then the cusp comes into the actual point which can be 
interpreted as the intersection of 
${\cal L}_{34}, {\cal L}_{35}$  and  ${\cal L}_{45}$.  
The boundary of an equidistant neighborhood of 
the point enjoys a spherical geometry.  
The global topology change from  
${\cal Q}(2/5,2/5,2/5,2/5,2/5)$  to  ${\cal Q}(\mu)$  can be 
described by collapsing down a $-1$ rational curve 
to a point, which is nothing but a blowing down.  
\end{Example} 

\begin{ThNum} 
Any topology change within a family 
${\cal Q}(\mu_1, \cdots, \mu_5)$  under the 
condition (\ref{Eq:weight})  is attained 
by a sequence of blowing up and down along  ${\cal L}_{ij}$'s. 
\end{ThNum}

\begin{pf} 
Since the total sum of the  $\mu_j$'s  equals  $2$,  
possible values of  $\mu_j$'s  such that some pair has 
the sum equal to  $1$  are limited.  
Either one pair does, two
pairs with a common value do or three values equal  $1/2$.  
In particular, at most three irreducible components 
of  ${\cal L}$, forced to be disjoint, 
are involved with cusp opening or closing.  
Hence the claim follows easily from this naive observation 
with Proposition \ref{Prop:ConeAngle}.  
\end{pf} 

\begin{LemNum} 
Suppose that  $1 \leq j < i \leq 5$  both are different 
from  $1 \leq k < l \leq 5$.  
When  $\mu_j + \mu_i > 1$,  then 
${\cal L}_{kl}$  is a hyperbolic cone sphere with 
cone singularity at blown down  ${\cal L}_{ji}$   
with cone angle $= - 2\pi(1-\mu_j-\mu_i)$.  
\end{LemNum}

\begin{pf} 
We may assume without loss of generality that  
$j = 1, i = 2, k = 3, l = 4$.  
Then  ${\cal L}_{34}$  can be identified with  
${\cal Q}(\mu_1, \mu_2, \mu_3+\mu_4, \mu_5)$  by definition.  
It has three cone points coming from the 
intersection with  
${\cal L}_{15}, \, {\cal L}_{25}$,  and 
${\cal L}_{35}$  and  ${\cal L}_{45}$  simultaneously 
which is appeared by blowing down  ${\cal L}_{12}$.  
Hence the cone angle around the last cone point is 
calculated as 
\begin{equation*} 
	2\pi(1-((\mu_3+\mu_4)+\mu_5)) 
	= 2 \pi((1-(2-\mu_1-\mu_2)) 
	= -2\pi(1-\mu_1-\mu_2).  
\end{equation*} 
by Proposition \ref{Prop:ConeAngle}.  
\end{pf} 

The following effectiveness of deformations 
should be compared with Example \ref{Ex:IsometricResult}.  

\begin{ThNum} 
Suppose  $n=5$  and the weight satisfies the 
condition (\ref{Eq:weight}).  
If  ${\cal Q}(\mu)$  is isometric to  ${\cal Q}(\lambda)$,  
then there is a permutation  $\sigma$  of five letters 
such that   $\sigma(\mu) = \lambda$.  
\end{ThNum}

\begin{pf}
Given the weight  $\mu$,  we get ten 
numerical invariants  $2\pi(1 - \mu_j - \mu_i)$  by 
running  $1 \leq j < i \leq 5$,  which describe 
cone angles appeared in  ${\cal Q}(\mu)$.  
If these ten numerical invariants are the same 
for  ${\cal Q}(\mu)$  and  ${\cal Q}(\lambda)$, 
then it is quite easy to check that 
the sets of components of  $\mu$  and  $\lambda$  must 
be the same.  
\end{pf}


\section{Problems}

Here we list a few problems arisen in the study. 

\begin{Problem} 
Work out a similar study in the last section for 
$n \geq 6$.  
\end{Problem} 

\begin{Problem} 
Develop a deformation theory of complex hyperbolic cone  
structures in complex dimension $2$ or higher, 
and show how much structures come from the weighted configurations.  
\end{Problem} 

The second problem derives a few subquestions.  
Since the complex hyperbolic geometry is rigid, 
the space we should look at is 
the space of representations with appropriate data.  
Suppose  $n = 5$  and recall 
Example \ref{Ex:TopologicallyStableDeformation}.  
Define the subspace 
\begin{equation*} 
	R \subset \Hom(\pi_1({\cal Q}), \PU(1,2))/{\text{conjugacy}}
\end{equation*} 
to be the set of representations up to conjugacy 
such that each meridional element of  
${\cal L}_{ji}$  in  $\pi_1({\cal Q})$  
is represented by appropriate rotational elements.  

\begin{Problem}
Dose  $R$  locally parameterize the deformations of 
complex hyperbolic cone structures base on  
$(X, {\cal L})$  in 
Example \ref{Ex:TopologicallyStableDeformation} ?  
If so, what is  $\dim R$  at the holonomy 
representation of  ${\cal Q}(\mu)$  under 
the condition (\ref{Eq:Stability}) ?  
\end{Problem}

The same question for the real slice  ${\cal Q}_{\real}$  of  
${\cal Q}$,  formed by 
the point configurations lying on the circle up to projective 
transformations, has been discussed in  
\cite{KojimaNishiYamashita, YamashitaNishiKojima}.  
When  $n = 5$,  ${\cal Q}_{\real}(\mu)$  is 
a nonsingular hyperbolic surface homeomorphic to 
a connected sum of  $5$  copies of the 
real projective surface  $\real \proj^2$  under 
the condition  (\ref{Eq:Stability}).  
The dimension of the set of deformations coming 
from the real slice of the weighted configurations 
is  $4$  though the dimension of the space 
of hyperbolic structures on such a surface is  $9$.  

\begin{Problem}
Is a complex hyperbolic cone structure based on 
a pair  
$(X, {\cal L})$  uniquely 
determined by the cone angles for  ${\cal L}_{ji}$'s ?
Moreover does the angle fixing rigidity hold 
for complex hyperbolic cone manifolds with 
cone angles $\leq 2\pi$ in general ?  
\end{Problem} 

The angle fixing local rigidity for real hyperbolic 
cone 3-manifolds with vertexless singularity 
such that cone angles all  $\leq 2\pi$  is 
proved in  \cite{HodgsonKerckhoff},  and 
a global rigidity for the same cone manifolds with 
cone angles all  $\leq \pi$  is proved in  \cite{Kojima}.  
Hence it is not quite wild to expect to have such rigidities, 
though the topological constrain and angle bound 
for the singularity should be taken into account.


\end{document}